\documentclass[12pt]{amsart}
\usepackage{amsbsy}
\begin{document}

\newtheorem{theorem}{Theorem}
\newtheorem{proposition}{Proposition}
\newtheorem{lemma}{Lemma}
\newtheorem{corollary}{Corollary}
\newtheorem{definition}{Definition}
\newtheorem{remark}{Remark}
\newcommand{\tex}{\textstyle}
\numberwithin{equation}{section} \numberwithin{theorem}{section}
\numberwithin{proposition}{section} \numberwithin{lemma}{section}
\numberwithin{corollary}{section}
\numberwithin{definition}{section} \numberwithin{remark}{section}
\newcommand{\ren}{\mathbb{R}^N}
\newcommand{\re}{\mathbb{R}}
\newcommand{\n}{\nabla}
\newcommand{\iy}{\infty}
\newcommand{\pa}{\partial}
\newcommand{\fp}{\noindent}
\newcommand{\ms}{\medskip\vskip-.1cm}
\newcommand{\mpb}{\medskip}
\newcommand{\AAA}{{\bf A}}
\newcommand{\BB}{{\bf B}}
\newcommand{\CC}{{\bf C}}
\newcommand{\DD}{{\bf D}}
\newcommand{\EE}{{\bf E}}
\newcommand{\FF}{{\bf F}}
\newcommand{\GG}{{\bf G}}
\newcommand{\oo}{{\mathbf \omega}}
\newcommand{\Am}{{\bf A}_{2m}}
\newcommand{\CCC}{{\mathbf  C}}
\newcommand{\II}{{\mathrm{Im}}\,}
\newcommand{\RR}{{\mathrm{Re}}\,}
\newcommand{\eee}{{\mathrm  e}}
\newcommand{\eb}{{\bf  e}}
\newcommand{\LL}{L^2_\rho(\ren)}
\newcommand{\LLL}{L^2_{\rho^*}(\ren)}
\renewcommand{\a}{\alpha}
\renewcommand{\b}{\beta}
\newcommand{\g}{\gamma}
\newcommand{\G}{\Gamma}
\renewcommand{\d}{\delta}
\newcommand{\D}{\Delta}
\newcommand{\e}{\varepsilon}
\newcommand{\var}{\varphi}
\newcommand{\lll}{\l}
\renewcommand{\l}{\lambda}
\renewcommand{\o}{\omega}
\renewcommand{\O}{\Omega}
\newcommand{\s}{\sigma}
\renewcommand{\t}{\tau}
\renewcommand{\th}{\theta}
\newcommand{\z}{\zeta}
\newcommand{\wx}{\widetilde x}
\newcommand{\wt}{\widetilde t}
\newcommand{\noi}{\noindent}
\newcommand{\uu}{{\bf u}}
\newcommand{\xx}{{\bf x}}
\newcommand{\yy}{{\bf y}}
\newcommand{\zz}{{\bf z}}
\newcommand{\aaa}{{\bf a}}
\newcommand{\cc}{{\bf c}}
\newcommand{\jj}{{\bf j}}
\newcommand{\ggg}{{\bf g}}
\newcommand{\UU}{{\bf U}}
\newcommand{\YY}{{\bf Y}}
\newcommand{\HH}{{\bf H}}
\newcommand{\GGG}{{\bf G}}
\newcommand{\VV}{{\bf V}}
\newcommand{\ww}{{\bf w}}
\newcommand{\vv}{{\bf v}}
\newcommand{\hh}{{\bf h}}
\newcommand{\di}{{\rm div}\,}
\newcommand{\ii}{{\rm i}\,}
\newcommand{\nn}{{\bf  n}}
\newcommand{\inA}{\quad \mbox{in} \quad \ren \times \re_+}
\newcommand{\inB}{\quad \mbox{in} \quad}
\newcommand{\inC}{\quad \mbox{in} \quad \re \times \re_+}
\newcommand{\inD}{\quad \mbox{in} \quad \re}
\newcommand{\forA}{\quad \mbox{for} \quad}
\newcommand{\whereA}{,\quad \mbox{where} \quad}
\newcommand{\asA}{\quad \mbox{as} \quad}
\newcommand{\andA}{\quad \mbox{and} \quad}
\newcommand{\withA}{,\quad \mbox{with} \quad}
\newcommand{\orA}{,\quad \mbox{or} \quad}
\newcommand{\atA}{\quad \mbox{at} \quad}
\newcommand{\onA}{\quad \mbox{on} \quad}
\newcommand{\ef}{\eqref}
\newcommand{\ssk}{\smallskip}
\newcommand{\LongA}{\quad \Longrightarrow \quad}
\def\com#1{\fbox{\parbox{6in}{\texttt{#1}}}}
\def\N{{\mathbb N}}
\def\A{{\cal A}}
\newcommand{\de}{\,d}
\newcommand{\eps}{\varepsilon}
\newcommand{\be}{\begin{equation}}
\newcommand{\ee}{\end{equation}}
\newcommand{\spt}{{\mbox spt}}
\newcommand{\ind}{{\mbox ind}}
\newcommand{\supp}{{\mbox supp}}
\newcommand{\dip}{\displaystyle}
\newcommand{\prt}{\partial}
\renewcommand{\theequation}{\thesection.\arabic{equation}}
\renewcommand{\baselinestretch}{1.1}
\newcommand{\Dm}{(-\D)^m}

\title
{\bf
Sturmian multiple zeros for Stokes and Navier--Stokes equations in
$\re^3$ via solenoidal Hermite polynomials}

\author {V.A.~Galaktionov}

\address{Department of Mathematical Sciences, University of Bath,
 Bath BA2 7AY, UK}
\email{vag@maths.bath.ac.uk}


\date{\today}


\begin{abstract}

The Cauchy problem for the 3D
{\em Stokes} and {\em Navier--Stokes equations},
  \be
  \label{0.1}
  \begin{matrix}
\uu_t =- \n p+\D \uu, \quad {\rm div} \, \uu=0 \inB \re^3\times
(-1,0], \quad \mbox{and} \ssk\\
 \uu_t+ (\uu \cdot \n)\uu=- \n p+\D \uu, \quad {\rm div} \, \uu=0
 \inB \re^3 \times (-1,0],
 \end{matrix}
  \ee
  where $\uu=[u,v,w]^T$ is the vector field  and $p$ is the
  pressure,
 is considered.
 A smooth bounded initial data $\uu_0(x)$, with ${\rm div}\, \uu_0=0$, are prescribed at
 $t=-1$.

The  problem of formation of {\em multiple zeros} at the point
$(x,t)=(0,0^-)$ of the components of $\uu(x,t)$ is considered. In
recent years, such a classic problem, which, for the 1D heat
equation, was solved by Sturm  in 1836, was under scrutiny for a
number of parabolic, hyperbolic, elliptic, and dispersion PDEs. As
usual, such an analysis
 gives insight into a ``microscopic blow-up scale
properties" of the equations \ef{0.1}
 under consideration. It is
shown that  formation of multiple zeros of solutions can follow
``self-focusing" of  nodal sets moving according to zero surfaces
of the corresponding {\em solenoidal Hermite polynomials} as
eigenfunctions of a rescaled adjoint Hermite operator. This is
always the case for the first  problem in \ef{0.1}, but not always
for the second one.
Using such blow-up asymptotics allows to state a {\em unique
continuation theorem}.

A similar phenomenon is studied for the
  well-posed Stokes and Burnett equations with the minus bi-Laplacian
 instead of the standard viscosity operator in \ef{0.1}.


\end{abstract}

\maketitle





\section{Introduction: towards micro-scale structure of smooth
solutions}
\label{S0}

\subsection{Stokes, Navier--Stokes, and Burnett equations}

We consider the Cauchy problem for the {\em three-dimensional}
{\em linear Stokes} and {\em Navier--Stokes equations} (NSEs)
  \be
  \label{1St}
\uu_t =- \n p+\D \uu, \quad {\rm div} \, \uu=0 \inB \re^3\times
(-1,0], \quad \mbox{and}
 \ee
 \be
 \label{2NSE}
 \uu_t+ (\uu \cdot \n)\uu=- \n p+\D \uu, \quad {\rm div} \, \uu=0
 \inB \re^3 \times (-1,0].
  \ee
  where $\uu=[u,v,w]^T$ is the vector field  and $p$ is the
  corresponding
  pressure.
 A smooth bounded initial data $\uu_0(x)$ are prescribed at
 $t=-1$, with ${\rm div} \, \uu_0=0$.

 The problem of formation of {\em multiple zeros} at the point
$(x,t)=(0,0^-)$ of the vector field $\uu(x,t)$ is considered. In
recent years, such a classic problem, which, for the 1D heat
equation, was solved by C.~Sturm (1836), was under scrutiny for a
number of parabolic, hyperbolic, elliptic, and dispersion PDEs;
see short surveys for each type of linear and nonlinear PDEs in
\cite{GalPet2m}, and also \cite{GMNSE}, containing a most recent
survey. As usual, such analysis
 gives insight into a ``microscopic scale
properties" of the equations \ef{0.1}
 under consideration.
 Indeed, such an microscopic blow-up approach  was and is key in
 classic regularity problems including fundamental questions of
 the regularity of
characteristic boundary points and singular points in potential
and other related PDE theory; we refer to  recent surveys in
\cite{GalPet2m, GalMazf(u), GMNSE}, where matching blow-up
techniques were principally used.

In the present paper, in a similar manner,
  we show that formation of multiple zeros of solutions can follow
``self-focusing" of  nodal sets moving according to zero surfaces
of the corresponding {\em solenoidal Hermite polynomials} as
eigenfunctions of a rescaled adjoint Hermite operator. It turns
out that this is always the case for the Stokes problem. For the
NSEs, the situation is shown to be more complicated, though
similar zeros do exist. In particular, this allows us to state a
non-standard {\em unique continuation theorem} for such problems.



For more clear expressing  our main ``blow-up" techniques and
their applicability in general  PDE theory, we develop in Appendix
A (the B one contains the corresponding Hermitian spectral
analysis) at the paper end,
 as a natural extension,
 a similar multiple zeros  analysis of the Cauchy problem for the
{\em fourth-order Stokes-like equations} and
 {\em well-posed Burnett
equations}:
 \be
 \label{1B}
  \begin{matrix}
 \uu_t =-\n p  - \D^2 \uu, \quad{\rm div} \, \uu=0 \inB
 \re^N \times (-1,0],
 \quad \uu(x,-1)=\uu_0(x),\ssk\ssk\ssk\\
 \uu_t+ (\uu \cdot \n)\uu=-\n p  - \D^2 \uu, \quad{\rm div} \, \uu=0 \inB
 \re^N \times (-1,0],
 \quad \uu(x,-1)=\uu_0(x),
  \end{matrix}
  \ee
   where initial data $\uu_0$ are sufficiently smooth and satisfy ${\rm div} \, \uu_0=0$.
Here, we have the {\em bi-harmonic} diffusion operators on the
right-hand side of the $\uu$-equations.
 It turns out that our general scheme describing multiple zeros
 analysis can be applied; however, requiring another
 non-self-adjoint spectral theory for the corresponding rescaled
 operator, where {\em generalized} solenoidal Hermite polynomials
 naturally occur (Appendix B).

\subsection{Leray blow-up rescaled variables: why Hermite polynomials occur on micro-scales}

For semilinear NSEs \ef{2NSE},
  we perform Leray's type \cite{Ler34}\footnote{In particular, Leray proposed
  not only to look at a
  self-similar blow-up
  as $t \to T^-$ but also at a further similarity extension for $t>T$, i.e., in the complementary limit
   $t \to T^+$, so that the blow-up factor $\sqrt{T-t}$ is replaced by $\sqrt{t-T}$.
See some   historical and further comments on Leray's blow-up
scenario of 1934
  for
  the 3D NSEs can be found in \cite[\S~2.2]{GalNSE}.}
   {\em nonstationary blow-up} scaling
   with the blow-up time $T=0$:
   \be
   \label{ll2}
   \tex{
   \uu(x,t)=\frac 1{\sqrt{-t}} \, \hat \uu(y,\t), \quad
 p(x,t)= \frac 1{(-t)}\, P(y,\t), \quad
    y = \frac x{\sqrt{-t}}, \quad \t=
   -\ln(-t) \to +\iy
 }
 \ee
 as $t \to 0^-$.
 This yields
 the
rescaled equations for $\hat \uu=(\hat u^1, \hat u^2, \hat u^3)^T$
and $P$,
 \be
  \label{S1}
   \tex{
    \hat \uu_\t= \D \hat \uu -
 \frac 12 \, (y \cdot \n) \hat \uu -\frac 12 \, \hat \uu-  (\hat \uu
\cdot \n)\hat \uu- \n P , \quad
  \di \hat \uu =0 \inB \re^3 \times \re_+.
  }
  \ee

\ssk

As a standard next step,  we exclude the pressure from the
equations \ef{S1},
 \be
  \label{ww1}
    \begin{matrix}
    \hat \uu_\t=\HH(\hat \uu) \equiv  \big(\BB^*- \frac 12 \, I \big)\hat \uu
 -  \mathbb{P}\,(\hat \uu
\cdot \n)\hat \uu
   \inB \re^3 \times \re_+, \ssk\ssk\\
   \mbox{where} \quad \mathbb{P} \vv=\vv - \n \D^{-1}(\n \cdot
   \vv) \quad(\| \mathbb{P} \|=1)\qquad\qquad\quad
  \end{matrix}
  \ee
  is the
  Leray--Hopf projector of $(L^2(\re^3))^3$ onto  the subspace
  $\{{\bf w} \in (L^2)^3: \,\, {\rm div}\, {\bf w}=0\}$ of solenoidal vector
   fields\footnote{Of course, using ${\mathbb P}$ in \ef{ww1}  emphasizes an
    unpleasant fact that the NSEs are
    a {\em nonlocal} parabolic problem, so that a somehow full use of
 order-preserving  properties of the semigroup is illusive; though some ``remnants" of the
  Maximum Principle
  for such {\em second-order} flows may remain and actually appear from time to
  time in some results
   (but these are completely illusive for more difficult fluid models \ef{1B}).}. Another representation is
   $
   \mathbb{P} \vv=(v_1-R_1 \s,v_2-R_2 \s,v_3-R_3 \s)^T,
   $
   where $R_j$ are the Riesz transforms, with symbols
   $\xi_j/|\xi|$, and $\s=R_1 v_1+R_2 v_2+R_3 v_3$.
  We then first apply $\mathbb{P}$
 to the original velocity equation in \ef{2NSE} and next
 use the blow-up rescaling \ef{ll2}. Using the fundamental
 solution of $\D$ in $\re^3$
  \be
  \label{FF55}
   \tex{
   b_3(y)= - \frac 1{4 \pi} \, \frac 1{|y|},
    }
    \ee
the operator in \ef{ww1} is written in the form of Leray's
formulation
\cite[p.~32]{Maj02} 
 \be
 \label{HHH21}
  \begin{matrix}
\HH(\hat \uu) \equiv  \big(\BB^*- \frac 12 \, I \big)\hat \uu -
(\hat \uu \cdot \n )\hat \uu + C_3 \int\limits_{\re^3} \frac
{y-z}{|y-z|^3}\,\, {\rm tr} (\n \hat \uu (z,\t))^2\, {\mathrm d}z,
\ssk \\
 \mbox{where} \quad
 {\rm tr} (\n \hat
\uu (z,\t))^2= \sum_{(i,j)} \,  \hat u_{z_j}^i \hat u_{z_i}^j
\andA C_3= \frac 1{4 \pi}.\quad
 \end{matrix}
 \ee

It follows from \ef{HHH21} that, in order to describe asymptotic
behaviour of small solutions near multiple zeros, as a first step,
one needs a spectral theory for the linearized {\em Hermite
operator} $\BB^*$ in a proper solenoidal functional space. Of
course, this belongs to classic theory of self-adjoint operators;
see Birman--Solomjak \cite[p.~48]{BS}. Moreover, in full capacity
and specially for the NSEs in $\re^3$ and $\re^2$, this theory was
already developed (though has been used in the opposite large-time
behaviour of solutions as $t \to +\iy$) with eigenfunctions of the
adjoint (in $L^2$-metric) operator $\BB$. However, since this case
is self-adjoint, we can directly use such a theory in the
complementary blow-up limit $t \to 0^-$.



\section{Hermitian  spectral theory of the linear rescaled
operator $\BB^*$:
 point spectrum and  solenoidal Hermite polynomials}
 \label{S2}


Thus, approaching the point $(0,0)$ in the blow-up manner
\ef{ll2}, one observes Hermite's operator $\BB^*$ as the principal
linear part of the rescaled equation \ef{S1}. Writing it in
 the corresponding divergent form,
 \be
 \label{mm1}
 \tex{
    \BB^* \vv
 \equiv \frac 1 {\rho^*}\, \n \cdot (\rho^* \n \vv ),
 }
 \ee
 where the weight is
$\rho^*(y)={\mathrm e}^{-\frac{|y|^2}4}>0$,
 we observe that the actual rescaled evolution is now restricted to the
 weighted $L^2$-space $L^2_{\rho^*}(\re^3)$, with the
 exponentially decaying
  weight $\rho^*(y)$.
 Here, $ \BB^*$ is the  (``adjoint") Hermite operator with
 the point spectrum \cite[p.~48]{BS}
  \be
  \label{bbb1}
   \tex{
   \s(\BB^*)= \big\{ \l_k= - \frac{k}2, \quad
   k=|\b|=0,1,2,...\big\} \quad (\mbox{$\b$ is a multiindex}),
    }
    \ee
    where each $\l_k$ has the multiplicity $\frac{(k+1)(k+2)}2$
     (for $N=3$).
     The corresponding complete and closed set of eigenfunctions
     $\Phi^*=\{\psi_\b^*(y)\}$
is composed from separable Hermite polynomials. Note another
important  property of Hermite polynomials:
 \be
 \label{bn100}
 \forall \, \psi_\b^*, \quad \mbox{any derivative} \,\, D^\g  \psi_\b^* \quad \mbox{is also
 an eigenfunction with $k=|\b|-|\g| \ge 0$}.
  \ee
 Recall that \cite{BS}
  \be
  \label{bn10}
 \mbox{polynomial set} \,\,\, \Phi^* \,\,\, \mbox{is complete and closed in \,
  $L^2_{\rho^*}(\re^3)$}.
 \ee

Further spectral properties are convenient to demonstrate using
the linear operator $\BB$,
 \be
 \label{BBB1}
  \tex{
  \BB= \D + \frac 12\, y \cdot \n + \frac 32\, I \inB
  L^2_\rho(\re^3) \whereA \rho = \frac 1{\rho^*},
  }
  \ee
  which is adjoint to $\BB^*$ in the dual $L^2$-metric. It has the same point spectrum and the
  corresponding eigenfunctions are multiple of the same Hermite
  polynomials according to the well-known {\em generating
  formula}:
   \be
   \label{BBB2}
    \tex{
    \psi_\b(y)= \frac {(-1)^{|\b|}} {\sqrt{\b !}}\, D^\b F(y)
    \equiv \psi_\b^*(y) F(y)
    \whereA F(y)= \frac 1{(4\pi)^{3/2}} \, {\mathrm e}^{-|y|^2/4}
    }
    \ee
    is the rescaled kernel of the fundamental solutions of
    $D_t-\D$ in $\re^3 \times \re_+$.
Then, the  bi-orthonormality  holds:
  \be
  \label{bn1}
  \langle \psi_\b^*, \psi_\g \rangle=\d_{\b \g} \quad \mbox{for
  any} \quad \b, \, \g,
   \ee
   where $\langle \cdot, \dot \rangle$ is the scalar product in
   $L^2(\re^3)$.
As is well known, this dual $L^2$-metric can be also treated as a
weighted one in $L^2_{\rho^*}(\re^3)$ (where $\BB^*$ becomes
symmetric):
 $$
 \tex{
 \langle \psi_\b^*, \psi_g \rangle \equiv \int_{\re^3} F(y) \,
 \psi_\b^* \psi_\g^* \, {\mathrm d}y \sim \langle \psi_\b^*,
 \psi_g^*
 \rangle_{\rho^*},
 }
 $$
 since $F(y) \sim \rho^*(y)$, up to a constant multiplier.
However, we prefer to keep the bi-orthonormality in the
non-symmetric form \ef{bn1}, since a similar condition occurs in
the principally non-symmetric Burnett cases, see \ef{Ortog}.



\ssk

Obviously, one needs to consider eigenfunction expansions in the
solenoidal restriction
 \be
 \label{bn11}
\hat L^2_{\rho^*}(\re^3)= L^2_{\rho^*}(\re^3)^3\cap\{ \di \vv=0\}.
 \ee
Indeed, among the polynomials $\Phi^*=\{\psi_\b^*\}$, there are
many that well-suit the solenoidal fields. Namely, introducing the
eigenspaces
 $$
 \Phi_k^*= {\rm Span}\,\{\psi_\b^*, \,\, |\b|=k\}, \quad k \ge 1,
  $$
in view of \ef{bn100}, $\di$ plays a role of a ``shift operator"
in the sense that
 \be
 \label{sh1}
 \di: \Phi^{*3}_k \to \Phi^*_{k-1}.
  \ee

We next define  the corresponding solenoidal eigenspaces as
follows:
 \be
 \label{Sol1}
  \tex{
  {\mathcal S}_k^*= \{\vv^*=[v_1^*,v_2^*,v_3^*]^T: \quad {\rm div} \, \vv^*=0,
  \,\,v_i^* \in \Phi_k^*\}\whereA {\rm dim}\,\, {\mathcal S}_k^*=k(k+2);
  }
  \ee
see \cite{Gal02, Gal02A, Gallay06} and further references therein.

The study \cite{Gal02} deals with {\em global} asymptotics as $t
\to +\iy$, where the adjoint operator $\BB$ given  in \ef{BBB1}
occurs. Since $\BB$ is self-adjoint in $L^2_\rho(\re^3)$, almost
all the results from \cite[Append.~A]{Gal02A} are applied to
$\BB^*$.
 For a full collection,
 see \cite{Brand09, Brand07} for further large-time  asymptotic expansions and self-similar solutions.
In particular, this made it possible to construct therein fast
decaying solutions of the NSEs on each 1D stable manifolds with
the asymptotic behaviour\footnote{We present here only the first
term of expansion; as usual in dynamical system theory, other
terms in the case of ``resonance" can contain $\ln t$-factors
({\em q.v.} \cite{Ang88} for a typical PDE application); this
phenomenon was shown to exist
 for the NSEs in $\re^2$ \cite[p.~236]{Gal02A}.}
 \be
 \label{mn1}
  \tex{
 \uu_\b(x,t) \sim  \, t^{\l_k-\frac 12} \, \vv_\b \big( \frac x{\sqrt
 t}\big)+... \asA t \to \iy, \,\,\,\mbox{where} \,\,\, \vv_\b= {\vv_\b^*}F \in
 {\mathcal S}_k
   }
   \ee
 are solenoidal eigenfunctions of $\BB$.
  Namely, taking
    \be
    \label{kl1}
    \tex{
    \vv=[v_1,v_2,v_3]^T \in {\mathcal S}_k, \,\,\, v_i \in {\Phi}_k= {\rm Span}\,\big\{\psi_\b=
    \frac{(-1)^{|\b|}}{\sqrt{\b !}}\, D^\b F(y), \, |\b|=k \big\},
 }
  \ee
  where $F$ stands for the rescaled Gaussian in \ef{BBB2}, 
we have that
 \be
 \label{kl2}
  \tex{
 \di \vv = (v_1)_{y_1} + (v_2)_{y_2}+ (v_3)_{y_3}= \di (\vv^* F)
 \equiv (\di \vv^*) F - \frac 12\, y \cdot \vv^* \, F.
  }
  \ee
 This establishes  a one-to-one correspondence between
 solenoidal eigenfunction classes ${\mathcal S}_k^*$ in \ef{Sol1}
 for $\BB^*$ and ${\mathcal S}_k$ in \ef{mn1}
 for $\BB$; see \ef{Hp0}--\ef{HP2} below
 for the first  eigenfunctions
 $\vv_\b= \vv^*_\b F$. Therefore,
  ${\rm dim}\,\, {\mathcal
   S}_k=k(k+2)$, etc.;
see details and rather involved  proofs of the asymptotics
\ef{mn1} for $k=1$ and 2 in \cite{Gal02}.

\ssk

In particular,  those solenoidal Hermite polynomial eigenfunctions
of $\BB^*$ can be chosen as follows \cite[p.~2166-69]{Gal02A}
(the choice is obviously not unique; normalization constants are
omitted):
 \be
 \label{Hp0} \underline{\l_0=0:} \quad
 \vv_0^*=[1,1,1]^T=\eb \quad (\mbox{the first solenoidal Hermite polynomial}),
  \ee
 \be
 \label{HP1}
  \tex{
   \underline{\l_1=- \frac 12:} \quad
  \vv_{11}^*=\left[\begin{matrix}0\\-y_3\\y_2 \end{matrix} \right], \quad
\vv_{12}^*=\left[\begin{matrix}y_3\\0\\-y_1 \end{matrix} \right],
\quad \vv_{13}^*=\left[\begin{matrix}-y_2\\y_1\\0 \end{matrix}
\right]\,\,({\rm dim}\, {\mathcal S}_1^*=3);
 }
 \ee
 \be
 \label{HP2}
  \begin{matrix}
   \underline{\l_2=- 1:} \,\,\,
  \vv_{21}^*=\left[\begin{matrix}4-y_2^2-y_3^2\\y_1y_2\\-y_1y_3 \end{matrix} \right],
  \,\,
\vv_{22}^*=\left[\begin{matrix}y_1y_2\\4-y_1^2-y_3^2\\-y_2y_3
\end{matrix} \right], \,\,
\vv_{23}^*=\left[\begin{matrix}y_1y_3\\-y_2y_3\\4-y_1^2-y_2^2
\end{matrix} \right],
 \ssk\ssk\\
  \vv_{24}^*=-\left[\begin{matrix}0\\-y_1y_3\\y_1y_2 \end{matrix} \right], \quad
\vv_{25}^*=-\left[\begin{matrix}y_2y_3\\0\\-y_2y_1 \end{matrix}
\right], \qquad\qquad\qquad  \quad 
 \ssk\ssk\\
  \vv_{26}^*=\left[\begin{matrix}-y_2 y_3\\y_2y_3\\y_1^2-y_2^2 \end{matrix} \right],
  \,\,
\vv_{27}^*=\left[\begin{matrix}y_1y_2\\y_3^2-y_1^2\\-y_2y_3
\end{matrix} \right], \,\,
\vv_{28}^*=\left[\begin{matrix}y_2^2-y_3^2\\-y_1y_2\\y_1y_3
\end{matrix} \right]\,\,\,({\rm dim}\, {\mathcal S}_2^*=8),  \quad \mbox{etc.}
  \end{matrix}
 \ee

We need the following final conclusion.
By \ef{bn10}, the set of vectors $\Phi^{*3}$ is complete and
closed\footnote{Note a standard result of functional analysis: all
reasonable  polynomials are  complete in any weighted $L^p$-space
with an exponentially decaying weight; see the analyticity
argument in Kolmogorov--Fomin \cite[p.~431]{KolF}.} in
$L^2_{\rho^*}(\re^3)^3$, so that
\be
 \label{bn12}
  \tex{
 \forall \, \vv \in L^2_{\rho^*}(\re^3)^3 \LongA \vv= \sum_{(\b)} c_\b
 \vv^*_\b, \quad \vv_\b^* \in \Phi^{*3}_k, \,\,\, k= |\b| \ge 0,
 }
  \ee
  where $c_\b$ are scalars and, in a natural way,  the multiindex $\b$ arranges summation over  all  solenoidal
  Hermite polynomials.
   In fact, the only vector expansion  coefficient in \ef{bn12} can be the first one, $\cc_0$, so, for convenience,
  we will may use the following vector notation:
  \be
  \label{cb1}
  \cc_0=[c_0^1,c_0^2,c_0^3]^T \LongA \cc_0
 \vv^*_0 \equiv
 [c_0^1 v_{0 1}^*,c_0^2v_{0 2}^*,c_0^3v_{0 3}^*]^T \quad (\vv_0^*=[1,1,1]^{T}).
 \ee


  Thus,  it then follows
  from \ef{bn1}--\ef{sh1} that
   \be
  \label{bn10N}
 \mbox{polynomial set} \,\,\, \hat \Phi^*= \Phi^{*3}\cap\{\di \vv=0\} \,\,\,
 \mbox{is complete and closed in \,
  $\hat L^2_{\rho^*}(\re^3)$}.
 \ee
In what follows, we always assume that we  deal with ``solenoidal"
asymptotics involving eigenfunctions as in \ef{Sol1}.



\ssk

For Burnett equations in \ef{1B}, as we have promised to go with
in parallel, the blow-up rescaling and elements of linear
solenoidal spectral theory are found in Appendices A and B.

\section{First application of Hermitian spectral theory: Sturmian local
structure of zero sets of bounded solutions and unique
continuation}
 \label{SHerm}

\subsection{A dynamical system for Fourier coefficients}



 Consider the NSEs \ef{2NSE}.
 We assume that, in a neighbourhood of the point $(x,t)=(0,0^-)$, the
solution $ \uu(x,t)$ is uniformly bounded and is such that the
eigenfunction expansion, as in \ef{bn12},
 \be
 \label{ex1}
 \tex{
  \hat \uu(y,\t)= \sum_{(\b)} c_\b(\t) \vv_\b^*(y),
  }
  \ee
 converges in $\hat L^2_{\rho^*}(\re^3)$, and moreover, uniformly
 on compact subsets. These convergence questions of polynomial series
  are  standard \cite{BS}; see also \cite{Eg4, 2mSturm}, where generalized Hermite polynomials
   occur and further references and  details are given.
 In particular,  if $\uu(x,t)$ remains bounded for all $t \in [-1,0]$, then,
 obviously, for such bounded data $\uu_0$,
  the
   convergence in \ef{ex1} always takes place.

 Then, the expansion coefficients satisfy the following dynamical
 system (DS):
  \be
  \label{ex2}
   \left\{
   \begin{matrix}
 \dot c_\b= \big(\l_\b-\frac 12 \big)c_\b
 + \sum_{(\a,\g)} d_{\a\g\b} c_\a c_\g
  \quad \mbox{for any} \,\,\, |\b|
\ge 0,
 \ssk\ssk\\
 \mbox{where} \quad
 d_{\a\g\b}=  -  \langle\mathbb{P}\,(\hat \vv_\a^* \cdot
\n)\hat \vv_\g^*, \vv_\b \rangle \quad \mbox{for all}\quad \a,\,
\g.\quad
 \end{matrix}
 \right.
 \ee
 Since \ef{ex1} is a standard eigenfunction expansion via Hermite polynomials of
 a given bounded smooth rescaled solution $\hat \uu(y,\t) \in H^2_{\rho^*}(\re^3)$,
the quadratic sum on the right-hand side converges.
  Recall
that, moreover, according to the blow-up scaling \ef{ll2}, we
actually deal with bounded and uniformly exponentially small
rescaled solutions satisfying
 \be
 \label{ssjj1}
 |\hat \uu(y,\t)| \le C \,{\mathrm e}^{- \frac \t 2} \quad \mbox{in} \quad \re^3 \times
 \re_+.
  \ee

The DS \ef{ex2} is difficult for a general study. For instance, it
contains the answer to the existence/nonexistence of the
$L^\iy$-blow-up  question (The Millennium Prize Problem,
\cite{Feff00}), i.e., whether there exists a Type-II blow-up at
the internal point $(0,0^-)$
(see a discussion
in \cite[\S~5]{GMNSE}).

\ssk

 For
regular points, the DS \ef{ex2} can provide us with a typical
classification of multiple zeros and nodal sets of solutions. Note
again that this kind of study was first performed by Sturm in 1836
for linear 1D parabolic equations \cite{St36}; see historical and
other details in \cite[Ch.~1]{2mSturm}.

\ssk

Thus, following these lines, we clarify local zero sets of
solutions of the NSEs at regular points.
  Assume that
 \be
 \label{ex4}
\uu(0,0)={\bf 0}.
 \ee
 In this connection,
 recall that the first eigenfunction of $\BB^*$ with $\l_\b=0$
 \be
 \label{ex3}
 \vv_0^*(y)  = [1,1,1]^T, 
 \ee
  is the only ones that have an empty nodal set.
 Then, bearing in mind the blow-up scaling term $(1-t)^{-\frac 12} \equiv {\mathrm e}^{\frac \t 2}$
 in \ef{ll2}, we have to assume that (here, we use the convention
 \ef{cb1})
  \be
  \label{ex5}
\cc_0(\t)={\bf 0}\,\,\,\mbox{or}\,\,\,   \cc_0(\t) \to {\bf 0}
\asA \t \to +\iy \,\,\, \mbox{exponentially faster than ${\mathrm
e}^{-\frac \t 2}$}.
    \ee



\subsection{Polynomial nodal sets for the Stokes equations}

 A first
clue to a correct understanding of the DS \ef{ex2} is given by the
Stokes equations \ef{1St},
 i.e., without the quadratic convection term.
   Then \ef{ex2} becomes linear diagonal and is easily solved:
 \be
 \label{FF2}
  \tex{
  \dot c_\b= \big(\l_\b-\frac 12 \big)c_\b
  \LongA
  c_\b(\t)=c_\b(0) {\mathrm e}^{- \frac {(1+|\b|) \t}2} \quad
  \mbox{for any}
  \quad |\b|\ge 0.
  }
  \ee
 Therefore, according to \ef{ex1} (and bearing in mind the completeness-closure
 of the Hermite polynomials),
  all possible multiple zero asymptotics for the Stokes problem (its local
``micro-scale turbulence") are described by finite solenoidal
Hermite polynomials, and the zero sets of rescaled velocity
components also asymptotically, as $\t \to +\iy$ (i.e., $t \to
0^-$) obey the nodal Hermite structures.



\subsection{Nodal sets for the Navier--Stokes equations}

 Consider the full nonlinear dynamical system \ef{ex2}, which on
 integration is
 \be
 \label{FF4}
  \tex{
   c_\b(\t)= c_\b(0) {\mathrm e}^{- \frac {(1+|\b|) \t}2}-{\mathrm e}^{- \frac {(1+|\b|) \t}2}
   \int\limits_0^\t \sum_{(\a,\g)} d_{\a\g\b} (c_\a c_\g)(s) {\mathrm e}^{\frac {(1+|\b|) s}2}
   \, {\mathrm d}s.
    }
    \ee
It follows that the nonlinear quadratic terms in \ef{FF4}, under
certain assumptions, can affect the rate of decay of solutions
near the multiple zero. As usual in calculus, the indeterminacies
in this integral quadratic term can be tackled by L'Hospital rule,
but this is technically is very difficult.

\ssk

Since we are mainly interested in the study of nodal structures of
solutions by using the eigenfunction expansion \ef{ex1}, we
naturally need to assume that it is possible to choose the leading
decaying term  (or a  linear combination of terms) in this sum  as
$\t \to +\iy$. Then obviously these leading terms will
asymptotically describe the Hermitian polynomial structure of
nodal sets as $t \to 0^-$.
 For PDEs with local nonlinearities, this is done in a
standard manner as in \cite[\S~4]{2mSturm}; in the nonlocal case,
this seems can cause technical difficulties.
 However, the DS \ef{ex2} looks (but illusionary) as being obtained
 from a problem with local nonlinearities. In other words, the
 nonlocal nature of the NSEs is hidden in \ef{ex2} in the
 structure of the quadratic sum coefficients $\{d_{\a\g\b}\}$,
 and this does not affect the nodal set behaviour for some classes of
 multiple zeros.
We will check this as follows:

\ssk

\noi\underline{\sc Resonance zeros}. We consider a ``resonance
class" of multiple zeros. Namely, let us assume there exist a
multiindex
 subset ${\mathcal B}$ and a function ${ h}(\t) \to { 0}$ such that
 \be
 \label{vv1}
  \begin{matrix}
  c_\b(\t) \sim { h}(\t) \asA \t \to +\iy
  \quad \mbox{for any $\b \in {\mathcal B}$}, \quad\ssk\ssk\\
|c_\b(\t)| \ll |{ h}(\t)| \asA \t \to +\iy
  \quad \mbox{for any $\b \not \in {\mathcal B}$}.
 \end{matrix}
 \ee
 In other words, only the coefficients $\{c_\b(\t), \, \b \in  {\mathcal
 B}\}$ are assumed to define the nodal set via \ef{ex1}, and other terms
 are negligible as $\t \to +\iy$.
 Under the natural assumption of a strong enough convergence of the quadratic sums
 in \ef{ex2} (this is expected not to be valid in the case of singular
 blow-up points only), taking the ODEs from \ef{ex2} for each $\b
 \in {\mathcal B}$ yields, for $\t \gg 1$,
  \be
  \label{hh1}
   \tex{
 \dot c_\b= \big(\l_\b-\frac 12 \big)c_\b + o(c_\b)
 \whereA c_\b(\t) \sim { h}(\t).
 }
 \ee
 Hence, the asymptotic balancing of these equations
must assume that as $\t \to + \iy$
 \be
  \label{hh2}
   \tex{
  \dot{ h} \sim \big(\l_\b-\frac 12 \big){ h}
 \LongA c_\b(\t) \sim { h}(\t) \sim {\mathrm e}^{(-\frac k2-\frac
 12)\t} \andA |\b|=k,
  }
  \ee
 where we may omit lower-order multipliers. Thus, there exists a $k \ge 1$ such that
  $|\b|=k$ for
 any $\b \in {\mathcal B}$.
 One can see that, for such ``resonance" multiple zeros, the
 nonlocal quadratic term in \ef{ex2} is not important.
 Thus, in the resonance zero class  prescribed by \ef{vv1}
as $\t \to +\iy$, on compact subsets in $y$, similar to Stokes'
problem,
 \be
 \label{ex6}
  \fbox{$
 \mbox{the nodal set of $\hat \uu(y,\t)$ is governed by some
 solenoidal {\bf Hermite} polynomials.}
  $}
  \ee

\ssk

\noi\underline{\sc Polynomial structure of multiple zeros is
universal}. Note that the conclusion that, locally, for any zero
of finite order at $(0,0)$,
  \be
  \label{ex111}
  \mbox{nodal sets of $ \uu(x,t)$ are governed by finite-degree polynomials}
   \ee
 is trivially true for any sufficiently smooth solution.
 Indeed, this follows from the Taylor expansion of such solutions
  \be
  \label{Tay44}
   \tex{
   \uu(x,t)= \sum\limits_{(|\mu|,|\nu| \le K)} C_{\mu\nu}\, x^\mu \, (-t)^\nu + {\mathbf R}_K(x,t)
   \whereA  C_{\mu\nu}=  \frac {(-1)^{\nu}}{\mu! \, \nu!}
   \big(D^{\mu,\nu}_{x,t}\uu\big)(0,0)
 }
  \ee
  and ${\mathbf R}_K=o(|x|^K(-t)^K)$ is a higher-order remainder.
  Translating \ef{Tay44} via \ef{ll2} into the expansion for $\hat
  \uu(y,\t)$ yields some polynomial structure, so \ef{ex111} is
  obviously true.
 Thus, the principal feature of \ef{ex6} is that the Hermite
 polynomials count only therein.

\ssk

\noi\underline{\sc Non-resonance zeros: a general classification}.
Obviously, for the nonlocal problem \ef{ww1}, there exist other
non-resonance zeros. Indeed, let $(0,0)$ be a zero of $\uu(x,t)$
of a finite order $M
\ge 1$, i.e.,  
 \be
 \label{mm111}
  \tex{
  \uu(x,0) = \sum_{(|\s|=M)} a_\s x^\s(1+o(1)) \sim x^\s \asA x \to
  0 \quad \big( \sum_{(|\s|=M)} |a_\s| \not = 0 \big).
  }
   \ee
We now use the following expansion:
 \be
 \label{mm2}
  \tex{
   \uu(x,t)= \uu(x,0)- \uu_t(x,0)(-t) + \frac 1{2!}\,
   \uu_{tt}(x,0)(-t)^2+...\,,
   }
   \ee
   where, by \ef{ww1}, all the time-derivatives
   $D_t^\mu \uu(x,0)$ can be calculated:
    \be
    \label{mm3}
 \uu_t(x,0)= \D \uu(x,0)+ (\mathbb{P}(\uu \cdot \n)\uu)(x,0)
 \sim x^{\s-2} + (\mathbb{P}(\uu \cdot \n)\uu)(x,0),
  \ee
  with a natural meaning of $\D x^\s \sim x^{\s-2}$.
 If the nonlocal term is negligible here and for other
 time-derivatives, i.e.,
  $$
 \uu_t(x,0) \sim x^{\s-2}, \quad  \uu_{tt}(x,0) \sim x^{\s-4}, ...
 \, ,
  $$
  then according to \ef{mm2} this leads to a Hermitian structure of
  nodal sets. In fact, this repeats
the  pioneering zero-set calculations performed  by Sturm (1836);
see his
 original computations in \cite[p.~3]{GalGeom}.

\ssk

In general,
 the nonlocal term in \ef{mm3} is not specified by a local
 structure of the zero under consideration, so, obviously, it can essentially
 affect the zero evolution. For instance, as a hint, we can have the
 following zero:
  \be
  \label{mm4}
   \tex{
   \uu_t(0,0)={\bf C} \not = 0\,\, \Rightarrow\,\,
    \uu(x,t) \sim \sum_{(|\s|=M)} a_\s x^\s - (-t) = {\mathrm
    e}^{-\t}\big(\sum_{(|\s|=M)} a_\s z^\s-1\big),
     }
     \ee
 where $z= \frac x{(-t)^{1/m}}$.
Hence, this nodal set is governed by the rescaled
 variable $z$, which is different from the standard similarity one
 $y$ in \ef{ll2}. Of course, due to the nonlocality of the
 equation, many other types of zeros can be described. Actually,
 such non-resonance zeros can be governed by sufficiently
 arbitrary finite polynomials as the general expansion \ef{Tay44}
 suggests. However, there exists a countable family of
 ``admissible" rescaled variables. Recalling that bounded smooth solutions of the NSEs are
 analytic in both $x$ and $t$ (see references below), we have to have that there exists a
 finite $K \ge 1$ such that
  \be
  \label{KK1}
  D_t^K \uu(0,0) \not = 0 \andA D_t^s \uu(0,0)=0 \forA
  s=1,2,...,K-1.
 \ee
Therefore, close to $(0,0)$, the structure of such a multiple zero
is given by
 \be
 \label{KK2}
  \tex{
  \uu(x,t) \sim
   \sum_{(|\s|=M)} a_\s x^\s - (-t)^K= {\mathrm e}^{-K \t}
  \big( \sum_{(|\s|=M)} a_\s z^\s-1 \big), \,\, z= \frac x{(-t)^\g}, \,\, \g=\frac K{M}.
  }
  \ee
Since $K$ and $M=|\s|$ are arbitrary positive integers, the
exponent $\g$ in the expansion \ef{KK2} can be an arbitrary
positive rational number.
 Thus, {\em the rescaled functions and
variables in $\ef{KK2}$ exhaust {\sc all} types of zero surfaces
(points) focusing as $x,\,t \to 0$  for the NSEs in $\re^3$}.


\ssk

  Finally, the proof that zeros of infinite order are
not possible for smooth non-analytic PDEs  (and, as usual in such
Carleman and Agmon-type uniqueness results, this occurs for $\uu
\equiv 0$ only) is a difficult technical problem; see an example
in \cite[\S~6.2]{2mSturm}. For analytic in $y$ solutions of the
NSEs
 (see  references and results in \cite{Dong07, Dong09, Zub07}),
 this problem is nonexistent, and then in \ef{ex6} the degree of
the solenoidal vector Hermite polynomials is always finite, though
can be arbitrarily large.

\subsection{An application: a unique continuation theorem}

Note another straightforward consequence of this analysis that
this gives the following conventional {\em unique continuations
result}: {\em let \ef{ex4} hold,  $(0,0)$ be a
resonance zero\footnote{Indeed, this is hard to check. However,
for the Stokes equations \ef{1St}, as well as for  any
sufficiently smooth PDEs with local nonlinearities (see
\cite{2mSturm}), this assumption is not needed, so that such a
unique continuation theorem makes a full sense.}, and
at least one component of the nodal set of $\hat \uu(y,\t)$ does
not obey \ef{ex6}. Then}
 \be
 \label{ex7}
  \uu \equiv {\bf 0} \quad \mbox{\em everywhere}.
   \ee
Of course, this is not that surprising, since the result is just
included in the existing and properly converging eigenfunction
expansion \ef{ex1} under the assumption \ef{vv1}.
According to \ef{Tay44}, there exists another ``funny version" of
the unique continuation result: \ef{ex7} {\em holds if a multiple
zero is formed in a non-polynomial self-focusing of zero surfaces,
or  via a rescaled variable not available in $\ef{KK2}$}, but this
is indeed trivial.

For elliptic equations $P(x,D)u=0$, this has the natural
counterpart on {\em strong unique continuation property} saying
that nontrivial solutions cannot have zeros of infinite order; a
result first proved by Carleman in 1939 for $P=-\D + V$, $V \in
L^\iy_{\rm loc}$, in $\re^2$ \cite{Carl39}; see \cite{Dos05,
Tao08} for further references and modern extensions.


\ssk

Thus, this is the first application of solenoidal Hermitian
polynomial vector fields for regular solutions of the NSEs.
 We expect that, due to the DS \ef{ex2}, some ``traces" of such an analysis and Hermite polynomials should be seen
in the fully nonlinear study of $\hat \uu(y,\t)$ at the singular
blow-up  point $(0,0)$, where, instead of \ef{ex4}, we have to
assume that, in the sense of $\limsup_{x,t}$,
 \be
 \label{ex8}
 |\uu(0,0)|=+\iy.
  \ee

\smallskip

\ssk

{\bf Acknowledgements.} The authors would like to thank
I.V.~Kamotski for useful discussions of various regularity issues
concerning  the NSEs and general PDE theory, as well as for his
advice to separate Sections 2.5, 2.6, and 3.1 from a
``discussion-survey" preprint \cite{GalNSE} to create the present
paper.




\begin{appendix}
\begin{small}
\section*{Appendix A: Multiple zeros for Burnett equations}
 \label{SBur}
 \setcounter{section}{1}
\setcounter{equation}{0}
\subsection{Burnett equations}

For both the systems \ef{1B}, the blow-up scaling \ef{ll2} is
replaced by
 \be
 \label{yyy}
  \tex{
   \uu(x,t)= (-t)^{-\frac 34}\, \hat \uu(y,\t), \quad
   y= \frac x{(-t)^{1/4}},
   }
   \ee
   so that
 the rescaled system \ef{ww1} takes a similar form
\be
  \label{ww1yy}
   \tex{
    \hat \uu_\t=\HH(\hat \uu) \equiv  \big(\BB^*- \frac 34 \, I \big)\hat \uu
 -  \mathbb{P}\,(\hat \uu
\cdot \n)\hat \uu
   \inB \re^3 \times \re_+.
 }
  \ee
 The spectral
 theory of the given here ``adjoint" operator
  \be
  \label{NBN1}
   \tex{
   \BB^*=- \D^2 - \frac 14\, y \cdot \n \whereA
   \s(\BB^*)=\big\{\l_\b=- \frac{|\b|}4, \, |\b| =0,1,2,...\big\}
 }
    \ee
    with eigenfunctions being generalized Hermite polynomials
is available in \cite{Eg4}; a solenoidal extension in the same
lines is needed. Necessary spectral theory of the operator pair
$\{\BB,\BB^*\}$ is developed below, in Appendix B.

 Therefore, under the same assumptions, the
polynomial structure of nodal sets is guaranteed for the
corresponding Stokes-like and  Burnett equations (resonance zeros)
in \ef{1B}; and, moreover,  for an arbitrary $2m$th-order
viscosity operator $-(-\D)^m \uu$ therein.

\ssk

\noi{\bf Remark: Burnett equations in a hierarchy of hydrodynamic
models.}
 The Burnett equations in \ef{1B} appear as the {\em second
approximation} (the NSEs \ef{2NSE} being the {\em first one}) of
the corresponding kinetic equations on the basis of Grad's method
in
 Chapman--Enskog expansions for hydrodynamics.
Namely,   Grad's method applied to kinetic
 equations,  by expanding the kernel of the integral operators involved
 via
  those with pointwise
 supports,
   yields, in addition to the classic operators of the Euler equations, other
  viscosity parts as follows:
 $$
  \mbox{$
  D_t \uu \equiv
  \uu_t +(\uu \cdot \n)\uu =  
  \sum\limits_{n=0}^\infty \e^{2n+1} \D^n(\mu_n \D \uu)+...= \e\big(\mu_0 \D \uu+ \e^2 \mu_1
  \D^2 \uu+...\big)+...\, ,
   $}
  $$
  where $\e>0$ is the {\em Knudsen number} Kn;  see details in
  Rosenau's  regularization approach,
  \cite{Ros89}. In a full model, truncating such series at $n=0$ leads to the
  Navier--Stokes equations (\ref{2NSE}) (with $\mu_0>0$), while $n=1$ is associated with the
    Burnett equations in \ef{1B}.

 Recall also that
Burnett-type equations, with a small parameter, appeared as
higher-order viscosity approximations of the Navier--Stokes
equations, represent  an effective tool for proving existence of
their weak (``turbulent" in Leray's sense) solutions; see Lions'
monograph \cite[\S~6, Ch.~1]{JLi}.
Note that the
 ``Problem
 on blow-up/non-blow-up for
Burnett equations in \ef{1B} at an inner point" starts from
dimensions $N = 7$; for $N \le 6$, there exists a unique global
smooth $L^2$-solution, \cite[\S~6]{G3KS}. It is expected  that
this open problem in $\re^7$ is not easier at all than the classic
{\em Millennium Prize One} for the NSEs in $\re^3$ \ef{2NSE}. In
both cases, a construction (or proving its nonexistence) of a
Type-II blow-up singularity is necessary, since, most plausibly, a
Type-I self-similar blow-up solutions are nonexistent,
\cite[App.~B]{GMNSE}.

\end{small}
\end{appendix}

\begin{appendix}
\begin{small}
\section*{Appendix B: Solenoidal  Hermitian spectral theory for $2m$th-order operators}
 \label{SAp2}
 \setcounter{section}{2}
\setcounter{equation}{0}

 We describe the necessary
 spectral properties of the linear $2m$th-order differential
 operator in $\ren$ ($m=2$ for the Burnett  equations in \ef{1B})
\begin{equation}
 \label{B1*}
  \tex{
 \BB^* = (-1)^{m+1} \D^{m}_y - \frac {1}{2m}
 \, y \cdot \n_y,
 }
 \end{equation}
and of its $L^2$-adjoint $\BB$
 given
by
  \begin{equation}
 \label{B1}
  \tex{
 \BB = (-1)^{m+1} \D_y^{m} + \frac {1}{2m}
 \, y \cdot \n_y +  \frac{N}{2m}\, I.
 }
 \end{equation}
 As we have seen, for $m=1$, \ef{B1*} and \ef{B1}  are classic Hermite self-adjoint  operators
 with completely known spectral properties, \cite[p.~48]{BS}.
 For any $m \ge 2$, both operators \ef{B1*} and \ef{B1},
 though looking very similar to those for $m=1$,
 {\em are not symmetric}  and do not admit a self-adjoint
extension, so we follow \cite{Eg4} in presenting  spectral theory.


\subsection{Fundamental solution, rescaled kernel, and first estimates}
 \label{Sect3}

  The fundamental solution $b(x,t)$  of the linear poly-harmonic parabolic equation
  \be
  \label{Lineq}
  u_t = - (-\D)^m u \inA
  \ee
 takes the standard similarity form
  \be
  \label{1.3R}
   \tex{
   b(x,t) = t^{-\frac N{2m}}F(y), \quad y= \frac  x{t^{1/{2m}}}.
    }
 \ee
 The rescaled kernel $F$ is the unique radial solution of the elliptic
 equation with the operator \ef{B1}, i.e.,
  \begin{equation}
\label{ODEf}
 {\bf B} F \equiv -(-\Delta )^m F + \textstyle{\frac 1{2m}}\, y \cdot
\nabla F + \textstyle{\frac N{2m}} \,F = 0
 \,\,\,\, {\rm in} \,\, \ren,  \quad \mbox{with} \,\,\, \textstyle{\int F =
 1.}
\end{equation}
 For  $m \ge 2$, the rescaled kernel function $F(|y|)$ is   oscillatory as $|y| \to \infty$ and
satisfies
 \cite{EidSys, Fedor}
\begin{equation}
\label{es11} 
 |F(y)| < D\,\,  {\mathrm e}^{-d_0|y|^{\alpha}}
\,\,\,{\rm in} \,\,\, \ren, \quad \mbox{where} \,\,\,
\a=\textstyle{ \frac {2m}{2m-1}} \in (1,2),
\end{equation}
for some positive constants $D$ and $d_0$ depending on $m$ and
$N$.

\subsection{Some constants}

 As we have seen,   the rescaled kernel $F(y)$
satisfies \ef{es11}, where $d_0$ admits an explicit expression;
see below.
 Such optimal exponential estimates of the fundamental solutions
 of higher-order parabolic equations are well-known and were first
 obtained by Evgrafov--Postnikov (1970) and Tintarev (1982); see
 Barbatis \cite{Barb, Barb04} for key references.

As a crucial  issue for the further boundary point regularity
study, we will need a sharper, than given by \ef{es11}, asymptotic
behaviour of the rescaled kernel $F(y)$ as $y \to +\iy$. To get
that, we keep four  leading terms in
\ef{ODEf} and
obtain, in terms of the radial variable $y \mapsto |y|>0$:
 \be
 \label{i1}
 \tex{
 (-1)^{m+1} \big[F^{(2m)} + m \frac{N-1}y \, F^{(2m-1)}+...\big]
  + \frac 1{2m} \, y F' + \frac N{2m} \, F=0 \forA y \gg 1.
 }
 \ee
 Using standard classic WKBJ asymptotics, we substitute into \ef{i1}
 the function
  \be
  \label{i2}
  F(y) = y^{-\d_0} \, {\mathrm e}^{a y^\a}+... \asA y \to + \iy,
   \ee
   exhibiting two scales.
   Balancing two leading terms
 gives the algebraic equation for $a$ and  $\d_0$:
 \be
 \label{i3}
  \tex{
 (-1)^m (\a a)^{2m-1}= \frac 1{2m} \andA \d_0=  \frac{m(2N-1)-N}{2m-1}>0\,
 .
 }
  \ee



By construction, one needs to get the root $a$ of \ef{i3} with the
maximal ${\rm Re}\, a<0$. This yields (see e.g., \cite{Barb,
Barb04})
 \be
 \label{i4}
  \tex{
 a= \frac{2m-1}{(2m)^\a} \big[-\sin\big( \frac{\pi}{2(2m-1)}\big) +
 \ii \cos\big( \frac{\pi}{2(2m-1)}\big)\big] \equiv -d_0 + \ii b_0
 \quad (d_0>0).
 }
 \ee
Finally, this gives the following double-scale asymptotic of the
kernel:
 \be
 \label{i5}
  \tex{
  F(y) =
   y^{-\d_0} \, {\mathrm e}^{-d_0 y^\a} \big[ C_1 \sin (b_0 y^\a)+
   C_2 \cos (b_0 y^\a)\big]+... \asA y=|y| \to + \iy ,
   }
   \ee
 where $C_{1,2}$ are real constants, $|C_1|+|C_2| \not = 0$.
 In \ef{i5}, we present the first two leading terms from the
 $m$-dimensional bundle of exponentially decaying asymptotics.

In particular, for the Burnett equations in \ef{1B} in $\re^3$, we
have
 \be
 \label{i6}
  \tex{
  m=2, \,\,N=3: \quad \a= \frac 43,  \quad d_0=3 \cdot 2^{-\frac{11}3},
 \quad b_0=3^{\frac 32} \cdot 2^{-\frac{11}3},
  \andA \d_0= \frac
  73.
   }
   \ee

\subsection{The discrete real spectrum and eigenfunctions of
 $\BB$}



 For $m \ge 2$,
 $\BB$ is considered  in the weighted space $L^2_\rho(\ren)$ with the
exponentially growing weight function
 \be
  \label{rho44}
  \rho(y) = {\mathrm e}^{a |y|^\a}>0 \quad {\rm in} \,\,\, \ren,
 \ee
  where $a \in (0,  2d_0)$ is a fixed
constant.
 We next
introduce a standard  Hilbert (a weighted Sobolev) space of
functions $H^{2m}_{\rho}(\ren)$ with the inner product and the
induced  norm
\[
 \tex{
 \langle v,w \rangle_{\rho} = \int\limits_{\ren} \rho(y) \sum\limits_{k=0}^{2m}
 D^{k}_y
 v(y) \, \overline {D^{k}_y w(y)} \,{\mathrm d} y \andA
\|v\|^2_{\rho} = \int\limits_{\ren} \rho(y) \sum\limits_{k=0}^{2m}
|D^{k}_y
 v(y)|^2 \, {\mathrm d} y.
  }
\]
Then $H^{2m}_{\rho}(\ren) \subset L^2_{\rho}(\ren) \subset
L^2(\ren)$, and  $\BB$ is a bounded linear operator from $
H^{2m}_{\rho}(\ren)$ to $ L^2_{\rho}(\ren)$. Key spectral
properties of the operator $\BB$ are as follows \cite{Eg4}:

\begin{lemma}
\label{lemspec}
 {\rm (i)}  The spectrum of $\BB$
comprises real simple eigenvalues only,
 \begin{equation}
\label{spec1}
 \tex{
 \sigma(\BB)=
\big\{\lambda_\b = -\frac k{2m}, \,\, k= |\b|= 0,1,2,...\big\}.
 }
\end{equation}
 {\rm (ii)} The eigenfunctions $\psi_\b(y)$ are given by
\begin{equation}
\label{eigen} \psi_\beta(y) =\textstyle{\frac{(-1)^{|\b|}}{\sqrt
{\b !}}} \, D^\beta F(y), \quad \mbox{for any} \,\,\,|\b|=k.
\end{equation}

\noi{\rm (iii)} Eigenfunction subset \ef{spec1} is complete
 in $L^2({\re})$ and in $L^2_{\rho}({\re})$.

\noi {\rm (iv)} The resolvent $(\BB-\lambda I)^{-1}$
for  $\lambda \not \in \sigma(\BB)$ is a compact integral operator
in $L^2_{\rho}(\ren)$.
\end{lemma}


By Lemma \ref{lemspec}, the   centre and stable subspaces of $\BB$
are given by
 \be
 \label{centr1}
E^c = {\rm Span}\{\psi_0= F\} \andA E^s = {\rm Span}\{\psi_\b, \,
|\b| \ge 1\}.
 \ee

\subsection{Polynomial eigenfunctions of the operator $\BB^*$}

Consider the operator (\ref{B1*}) in the weighted space
$L^2_{\rho^*}(\ren)$, where $\langle \cdot, \cdot
\rangle_{\rho^*}$ and $\|\cdot\|_{\rho^*}$ are
 the inner product
and the norm,
 with the   ``adjoint"  exponentially decaying weight
function
  \begin{equation}
\label{rho2}
 \tex{
 \rho^*(y) \equiv \frac  1 {\rho(y)} = {\mathrm e}^{-a|y|^{\a}} > 0.
 }
\end{equation}
 We ascribe to $\BB^*$ the domain
 $H^{2m}_{\rho^*}(\ren)$, which is dense in $L^2_{\rho^*}(\ren)$, and
 then
$$
 \BB^*: \,\, H^{2m}_{\rho^*}(\ren) \to L^2_{\rho^*}(\ren)
$$
 is a bounded linear operator.  $\BB$ is adjoint
 to $\BB^*$ in the usual sense: denoting by $\langle \cdot,\cdot \rangle $  the
inner product in the dual space $L^2(\ren)$, we have
  \begin{equation}
 \label{Badj1}
 \langle \BB v, w \rangle =  \langle v, \BB^* w \rangle
 \quad \mbox{for any} \,\,\, v \in H^{2m}_\rho(\ren) \andA
 w \in H^{2m}_{\rho^*}(\ren).
 \end{equation}
The eigenfunctions of $\BB^*$ take a particularly simple finite
polynomial form and are as follows:


 \begin{lemma}
\label{lemSpec2}
 {\rm (i)} $ \sigma(\BB^*)=\s(\BB)$.

 \noi{\rm (ii)} The eigenfunctions  $\psi^*_\b(y)$ of $\BB^*$ are generalized Hermite
 polynomials of degree $|\b|$ given by
 \begin{equation}
 \label{psi**1}
  \tex{
 \psi_\b^*(y) = \frac 1{\sqrt{\beta !}}
 \Big[ y^\b + \sum_{j=1}^{[|\b|/2m]} \frac 1{j !}(-\Delta)^{m j} y^\b
 \Big] \quad \mbox{for any} \quad \b.
 }
 \end{equation}

 \noi{\rm (iii)} Eigenfunction subset \ef{psi**1} is complete   in $L^2_{\rho^*}(\ren)$.

  \noi {\rm (iv)}
$\BB^*$ has a compact resolvent $(\BB^*-\lambda I)^{-1}$ in
$L^2_{\rho^*}(\ren)$ for $\lambda \not \in \sigma(\BB^*)$.

\noi{\rm (v)}  The bi-orthonormality of the bases $\{\psi_\b\}$
and $\{\psi_\g^*\}$ holds in the dual $L^2$-metric:
 \be
 \label{Ortog}
\langle \psi_\b, \psi_\g^* \rangle = \d_{\b\g} \quad \mbox{for
any} \quad \b,\, \g.
 \ee
\end{lemma}





\ssk

\noi{\bf Remark on closure.} This is an important issue for using
eigenfunction expansions of solutions.
 Firstly,
  in the self-adjoint case $m=1$, the sets of eigenfunctions are closed in the
corresponding spaces, \cite{BS} (and we have used this in our
previous NSEs study).

Secondly, for $m \ge 2$, one needs some extra details. Namely,
using (\ref{Ortog}), we can introduce the subspaces of
eigenfunction expansions and begin with the operator $\BB$. We
denote by $\tilde L^2_\rho$ the subspace of eigenfunction
expansions $v= \sum c_\b \psi_\b$ with coefficients $c_\b =
\langle v, \psi^* \rangle$ defined as the closure of the finite
sums $\{\sum_{|\b| \le M} c_\b \psi_\b\}$ in the norm of
$L^2_\rho$. Similarly, for the adjoint operator $\BB^*$, we define
the subspace $\tilde L^2_{\rho^*} \subseteq L^2_{\rho^*}$. Note
that since the operators are not self-adjoint and the
eigenfunction subsets are not orthonormal, in general, these
subspaces can be different from
$
 L^2_{\rho}$ and $L^2_{\rho^*}$, and particularly
the equality is guaranteed in the self-adjoint case $m=1$,
$a=\frac 1 4$.

Thus, for $m \ge 2$, in the above subspaces obtained via a
suitable closure, we can apply standard eigenfunction expansion
techniques as in the classic self-adjoint case $m=1$.


%

\subsection{Solenoidal Hermite polynomials}

The vector solenoidal Hermite polynomials are constructed from
\ef{psi**1} in a manner similar to that for $m=1$; cf
\ef{Hp0}--\ef{HP2}. Namely, given a vector polynomial
 \be
 \label{t6}
 \vv_\b^*=[\psi_{\b_1}^*, \psi_{\b_2}^*,..., \psi_{\b_N}^*]^T
 \whereA |\b_1|=|\b_2|=...=|\b_N|=|\b|,
  \ee
 it gets solenoidal provided that
  \be
  \label{t61}
   \tex{
   {\rm div}\, \vv_\b^* \equiv \sum\limits_{i=1}^N (\psi_{\b_i})_{y_i}
   =0.
   }
   \ee

For instance, for the Burnett case $m=2$ and $N=3$ (not all
linearly independent ones are presented, normalization constants
are omitted):
 \be
 \label{t62}
 \l_0=0: \quad \vv_0^*=[1,1,1]^T,
  \ee
  \be
  \label{t63}
   \tex{
   \l_1=-\frac 14: \quad \vv_{11}^*=[y_2, -y_3, y_2]^T, \quad
   \vv_{12}^*=[y_3,y_3,-y_1]^T, \quad \vv_{13}^*=[-y_2,y_1,y_1]^T,
   }
    \ee
 \be
 \label{t64}
  \tex{
 \l_2=- \frac 12: \quad \vv_{21}^*=[-y_1^2-y_3^2,y_1y_2,y_1y_3]^T, \,\,
 \vv_{22}^*=[y_1y_2,-y_2^2-y_3^2,y_2y_3]^T, \,\,\, \mbox{etc.}
 }
  \ee
  \be \label{t65}
   \tex{
  \l_3=- \frac 34: \quad \vv_{31}^*=[y_2^3,y_3^3,y_1^3], \quad
   \vv_{32}^*=
  [y_1y_2^2, y_2y_1^2,-y_3(y_1^2+y_2^2)], \,\,\,\mbox{etc.}
  }
   \ee
    \be
    \label{t66}
    \l_4=-1: \quad \vv_{41}^*=[y_2^4+4!,y_3^4+4!,y_1^4+4!]^T, \,\,\,
     \vv_{42}^*=[y_1y_2^3,y_2y_1^3,-y_3(y_1^3+y_2^3)]^T, \quad
     \mbox{etc.}
     \ee
As in the self-adjoint case $m=1$, some technical efforts are
necessary towards   completeness/closure of  generalized
solenoidal Hermite polynomials in suitable spaces. We omit
details.
\end{small}
\end{appendix}

\end{document}